\documentclass{amsart}%
\usepackage{amssymb}
\usepackage{amsfonts}
\usepackage{amsmath}
\makeatother
\newtheorem{theorem}{Theorem}

\newtheorem{corollary}[theorem]{Corollary}

\newtheorem{definition}[theorem]{Definition}

\newtheorem{lemma}[theorem]{Lemma}

\newtheorem{proposition}[theorem]{Proposition}
\newtheorem{remark}[theorem]{Remark}

\begin{document}

\title[Stone spectra of von Neumann algebras of type $I_{n}$]{Stone spectra of von Neumann algebras\\of type $I_{n}$}
\author[Andreas D\"{o}ring]{Andreas D\"{o}ring\\IAMPh, Fachbereich Mathematik\\J. W. Goethe-Universit\"{a}t Frankfurt, Germany}
\email{adoering@math.uni-frankfurt.de}
\date{17. January 2005}

\begin{abstract}
The Stone spectrum of a von Neumann algebra is a generalization of the Gelfand spectrum, as was shown by de Groote. In this article we clarify the structure of the Stone spectra of von Neumann algebras of type $I_{n}$.
\end{abstract}

\maketitle

\section{Introduction}

In a new approach tying together classical and quantum observables, de Groote
has developed the theory of \emph{Stone spectra} of lattices and
\emph{observable functions}. In this section, we will give the relevant
definitions and cite some of de Groote's results. For details, the reader is
referred to de Groote's work \cite{deG01,deG05}.

\ 

When speaking of a lattice, we will always mean a $\sigma$-complete lattice at
least. A lattice $\mathbb{L}$ always has a zero element $0$ and a unit element
$1$. The starting point is the new notion of a \emph{quasipoint} of a lattice,
which is nothing but a maximal filter base, generalizing to arbitrary lattices
what Stone did in the 1930s for Boolean algebras \cite{Sto36}:

\begin{definition}
\label{DQPInLat}A subset $\mathfrak{B}$ of a lattice $\mathbb{L}$ is called a
\textbf{quasipoint%
\index{quasipoint} of }$\mathbb{L}$%
\index{quasipoint!of a lattice} if it has the following properties:

\begin{enumerate}
\item[(i)] $0\notin\mathfrak{B},$

\item[(ii)] $\forall a,b\in\mathfrak{B}\ \exists c\in\mathfrak{B}:c\leq
a\wedge b,$

\item[(iii)] $\mathfrak{B}$ is maximal with respect to (i) and (ii).
\end{enumerate}
\end{definition}

It is easily seen that for a quasipoint $\mathfrak{B}$ of the lattice
$\mathbb{L}$, we have%
\[
\forall a\in\mathfrak{B}\ \forall b\in\mathbb{L}:(a\leq b\Longrightarrow
b\in\mathfrak{B}).
\]
In particular, $\forall a,b\in\mathfrak{B}:a\wedge b\in\mathfrak{B}$, so
quasipoints are maximal dual ideals also \cite{Bir73}. The set of quasipoints
of a lattice $\mathbb{L}$ is denoted by $\mathcal{Q}(\mathbb{L})$ and is
equipped with a natural topology (also inspired by Stone): for $a\in
\mathbb{L}$, let%
\[
\mathcal{Q}_{a}(\mathbb{L}):=\{\mathfrak{B}\in\mathcal{Q}(\mathbb{L}%
)\ |\ a\in\mathfrak{B}\}.
\]
Obviously, we have%
\[
\mathcal{Q}_{a\wedge b}(\mathbb{L})=\mathcal{Q}_{a}(\mathbb{L})\cap
\mathcal{Q}_{b}(\mathbb{L}),
\]
so $\{\mathcal{Q}_{a}(\mathbb{L})\ |\ a\in\mathbb{L}\}$ is the base of a
topology on $\mathbb{L}$.

\begin{definition}
The set $\mathcal{Q}(\mathbb{L})$ of quasipoints of a lattice $\mathbb{L}$,
equipped with the topology given by the sets $\mathcal{Q}_{a}(\mathbb{L)}$
defined above, is called the \textbf{Stone spectrum%
\index{Stone spectrum} of the lattice }$\mathbb{L}$%
\index{Stone spectrum!of a lattice}.
\end{definition}

One can show that the Stone spectrum $\mathcal{Q}(\mathbb{L})$ is a
zero-dimensional, completely regular Hausdorff space. For the example
$\mathbb{L}=\mathbb{L}(\mathcal{H})$, the lattice of closed subspaces of a
complex separable Hilbert space $\mathcal{H}$, $\mathcal{Q}(\mathbb{L})$ is
not compact if $\dim\mathcal{H}>1$. If $\mathcal{H}$ is infinite-dimensional,
then $\mathcal{Q(}\mathbb{L)}$ is not even locally compact \cite{deG01}. The
lattice $\mathbb{L}(\mathcal{H})$ plays an important role in the foundations
of quantum theory, which was first recognized by Birkhoff and von Neumann
\cite{BirvNeu36}. $\mathbb{L}(\mathcal{H})$ is isomorphic the $\mathcal{P(H)}%
$, the lattice of projections onto closed subspaces of $\mathcal{H}$.

\ 

Let $\mathcal{R}$ be a unital von Neumann algebra, given as a subalgebra of
the algebra $\mathcal{L(H)}$ of bounded operators on some Hilbert space
$\mathcal{H}$. The Stone spectrum $\mathcal{Q(R)}$ of $\mathcal{R}$ means the
Stone spectrum $\mathcal{Q(P(R))}$ of the projection lattice $\mathcal{P(R)}$
of $\mathcal{R}$. In quantum theory, including quantum mechanics in the von
Neumann representation, quantum field theory and quantum information theory,
\emph{observables%
\index{observable}%
\index{observable!physical}} are represented by self-adjoint operators $A$ in
some von Neumann algebra $\mathcal{R}$, the \emph{algebra of observables%
\index{observable!algebra of -s}}. The set of observables $\mathcal{R}_{sa}$
forms a real linear space in the algebra $\mathcal{R}$.

\ 

De Groote shows that if $\mathcal{R}$ is abelian, there is a homeomorphism
between the Stone spectrum $\mathcal{Q(R)}=\mathcal{Q(P(R))}$ and the Gelfand
spectrum $\Omega(\mathcal{R})$ of $\mathcal{R}$. Hence, for an arbitrary von
Neumann algebra $\mathcal{R}$, the Stone spectrum $\mathcal{Q(R)}$ is a
\emph{generalization of the Gelfand spectrum}.

\ 

Observable functions are introduced in the following way:

\begin{definition}
\label{DOFOfSelfAdjOp}Let $A\in\mathcal{R}_{sa}$, and let $E^{A}=(E_{\lambda
}^{A})_{\lambda\in\mathbb{R}}$ be the spectral family of $A$. The function%
\[
f_{A}:\mathcal{Q(R)}\longrightarrow\mathbb{R}%
\]
defined by%
\[
f_{A}(\mathfrak{B}):=\inf\{\lambda\in\mathbb{R}\ |\ E_{\lambda}^{A}%
\in\mathfrak{B}\}
\]
is called the \textbf{observable function corresponding to }$A$%
\index{observable function}%
\index{observable function!of self-adjoint operator}. $\mathcal{O}%
b(\mathcal{R}):=\{f_{A}\ |\ A\in\mathcal{R}_{sa}\}$ denotes the set of
observable functions of $\mathcal{R}$.
\end{definition}

One can show that the image of $f_{A}$ is the spectrum of $A$. Observable
functions are continuous functions, so%
\[
\mathcal{O}b(\mathcal{R})\subseteq C_{b}(\mathcal{Q(R)},\mathbb{R}),
\]
where $C_{b}(\mathcal{Q(R)},\mathbb{R})$ denotes the set of continuous bounded
real-valued functions on the Stone spectrum $\mathcal{Q(R)}$ of $\mathcal{R}$.
De Groote shows that equality only holds for abelian von Neumann algebras
$\mathcal{R}$. Moreover, if $\mathcal{R}$ is abelian, then the observable
functions turn out to be the Gelfand transforms of the self-adjoint elements
of $\mathcal{R}$, so for an arbitrary von Neumann algebra, the observable
functions are \emph{generalized Gelfand transforms}.

\ 

There are many more results on Stone spectra and observable functions and
their relation to physics \cite{deG01,deG05}. For example, observable
functions can be characterized intrinsically, without reference to
self-adjoint operators. Stone spectra also play some role in the proof of the
generalized Kochen-Specker theorem, which is an important no-go theorem on
hidden variables in quantum theory \cite{KocSpe67,Doe04}.

\ 

The elements of the Stone spectrum $\mathcal{Q(R)}$, the quasipoints, are
defined using Zorn's lemma. As usual, it is not easy to get some intuition of
such objects. Some extra structure of the von Neumann algebra is needed to
clarify the properties of the quasipoints and the Stone spectrum. The
quasipoints and the Stone spectrum of type-$I_{n}$-\emph{factors}
$\mathcal{R}=\mathcal{L(H)}$, where $\dim\mathcal{H}=n\in\{0,1,...\}$ is
finite, are known \cite{deG01}. We need the following: an isolated point
$\mathfrak{B}$ of the Stone spectrum $\mathcal{Q(R)}$ is called an
\emph{atomic quasipoint}. If $\mathcal{R}=\mathcal{L(H)}$, then the atomic
quasipoints of $\mathcal{P(R)}$($=\mathcal{P(H)}$) are of the form%
\[
\mathfrak{B}_{\mathbb{C}x}=\{P\in\mathcal{P(H)}\ |\ P\geq P_{\mathbb{C}x}\},
\]
where $x\in\mathcal{H}$, $|x|=1$ and $P_{\mathbb{C}x}$ is the projection onto
the line $\mathbb{C}x$. While for infinite-dimensional $\mathcal{H}$,
$\mathcal{P(H)}$ also has non-atomic quasipoints, for finite-dimensional
$\mathcal{H}$ the situation is simple:

\begin{proposition}
\label{PQPsOfTypeI_nFactors}If $\mathcal{H}$ is finite-dimensional, there are
only atomic quasipoints in $\mathcal{P(H)}$, and we have%
\[
\mathcal{Q(P(H))}=\{\mathfrak{B}_{\mathbb{C}x}\ |\ x\in S^{1}(\mathcal{H)}\},
\]
where $S^{1}(\mathcal{H)}$ denotes the unit sphere in Hilbert space. For
finite-dimensional $\mathcal{H}$, $\mathcal{P(H)}\simeq\mathbb{L}%
(\mathcal{H})$ is the projection lattice $\mathcal{P(R)}$ of a type $I_{n}$
factor $\mathcal{R}$, where $n=\dim\mathcal{H}$. This type $I_{n}$ factor
simply is represented as $\mathbb{M}_{n}(\mathbb{C})$, the $n\times n$ complex
matrices acting on $\mathcal{H}$.

Of course, the inner product of $\mathcal{H}$ plays no role here, but only the
linear structure, so we have also characterized the quasipoints of \ the
lattice of subspaces of a finite-dimensional vector space.
\end{proposition}

There were no results on the structure of quasipoints and Stone spectra of
more general von Neumann algebras up to now. In this article, we will examine
von Neumann algebras of type $I_{n}$ for finite $n$. Type $I_{n}$ algebras
include all von Neumann algebras on finite-dimensional Hilbert spaces and all
abelian von Neumann algebras. (The latter are of those of type $I_{1}$.) We
will make use of the fact that such an algebra is of the form $\mathbb{M}%
_{n}(\mathcal{A})$, where $\mathcal{A}$ is the center of $\mathcal{R}$. Note
that while $\mathcal{R}\simeq\mathbb{M}_{n}(\mathcal{A})$ is given by
\textquotedblleft finite\textquotedblright\ $n\times n$-matrices, the center
$\mathcal{A}$ of $\mathcal{R}$ may be represented on an infinite-dimensional
Hilbert space.

\ 

Drawing on a result on \emph{abelian quasipoints} (section \ref{_AbQPsOfVNAs}%
), i.e. quasipoints containing an abelian projection, a fairly complete
characterization of the Stone spectrum of a type $I_{n}$ algebra is obtained.
In particular, all quasipoints of a type $I_{n}$ algebra are abelian (Thm.
\ref{TTypeI_nAllQPab}) and the orbits of the action of the unitary group on
the Stone spectrum $\mathcal{Q(R)}$ are parametrized by the quasipoints of the
center of $\mathcal{R}$ (Thm. \ref{TQPsOfCenterParamUnitaryOrbits}).

\section{Abelian quasipoints of von Neumann algebras\label{_AbQPsOfVNAs}}

In this section, we will regard quasipoints containing an abelian projection.
It will be shown that there is a close relationship between the abelian
quasipoints of a von Neumann algebra $\mathcal{R}$ and the quasipoints of the
center of $\mathcal{R}$. This result will be central to the classification of
Stone spectra of type $I_{n}$ von Neumann algebras.

\begin{definition}
\label{DAbQP}A quasipoint $\mathfrak{B}\subseteq\mathcal{P(R)}$ is called
\textbf{abelian%
\index{quasipoint!abelian}} if it contains an abelian projection
$E\in\mathcal{R}$. The set of abelian quasipoints of a von Neumann algebra
$\mathcal{R}$ is denoted by $\mathcal{Q}^{ab}\mathcal{(R)}$.
\end{definition}

\begin{definition}
\label{DETrunk}The $E$\textbf{-trunk%
\index{quasipoint!e@$E$-trunk of}} $\mathfrak{B}_{E}$ ($E\in\mathfrak{B}$) of
a quasipoint $\mathfrak{B}$ is the set%
\[
\mathfrak{B}_{E}:=\{F\in\mathfrak{B}\ |\ F\leq E\}.
\]
Obviously, $\mathfrak{B}_{E}$ is a filter basis.
\end{definition}

\begin{lemma}
\label{LETrunkDeterminesQP}The $E$-trunk $\mathfrak{B}_{E}$ uniquely
determines the quasipoint $\mathfrak{B}$.
\end{lemma}

\begin{proof}
Let $\mathfrak{B}_{1},\mathfrak{B}_{2}$ be two quasipoints whose $E$-trunk is
$\mathfrak{B}_{E}$. Let $F$ be a projection in $\mathfrak{B}_{1}$. Then we
have $E\wedge F\in\mathfrak{B}_{E}\subset\mathfrak{B}_{2}$. If a quasipoint
contains a projection, it contains all larger projections, so $F\in
\mathfrak{B}_{2}$ and $\mathfrak{B}_{1}=\mathfrak{B}_{2}$ follows.
\end{proof}

\ 

This lemma holds analogously for any lattice $\mathbb{L}$, since no features
of the von Neumann algebra are used.

\begin{definition}
\label{DPartIsomOnQP}Let $\mathcal{R}\subseteq\mathcal{L(H)}$ be a von Neumann
algebra, $\mathfrak{B}\subset\mathcal{Q}_{E}(\mathcal{R})$ a quasipoint
containing $E$ and $\theta\in\mathcal{R}$ a partial isometry such that
$E=\theta^{\ast}\theta$. We set%
\[
\theta(\mathfrak{B}_{E}):=\{\theta F\theta^{\ast}\ |\ F\in\mathfrak{B}_{E}\}.
\]

\end{definition}

\begin{lemma}
If $\mathcal{R}\subseteq\mathcal{L(H)}$ is a von Neumann algebra and
$\theta\in\mathcal{R}$ is a partial isometry such that $E:=\theta^{\ast}%
\theta$, then for all projections $P_{U}\in\mathcal{R}$ such that $P_{U}\leq
E$ it holds that%
\[
\theta P_{U}\theta^{\ast}=P_{\theta U}.
\]

\end{lemma}

\begin{proof}
For $x\in U$, we have%
\[
\theta P_{U}\theta^{\ast}\theta x=\theta P_{U}Ex=\theta x=P_{\theta U}\theta
x.
\]
If $y\in(\theta U)^{\perp}$, then $\theta^{\ast}y\in U^{\perp}$ and thus%
\[
\theta P_{U}\theta^{\ast}y=0=P_{\theta U}y.
\]

\end{proof}

\ 

Obviously, $\theta(\mathfrak{B}_{E})$ is the $\theta E\theta^{\ast}$-trunk of
a quasipoint of $\mathcal{R}$ (given that $\theta^{\ast}\theta=E$). We will
denote the quasipoint induced by $\theta E\theta^{\ast}$ by $\theta
_{\mathcal{Q}}(\mathfrak{B}_{E})$.

\begin{remark}
Since in general $\theta^{\ast}\theta\notin\mathfrak{B}$ for an arbitrary
partial isometry $\theta$ and an arbitrary quasipoint $\mathfrak{B}$, we have
no action of the set of partial isometries on the Stone spectrum
$\mathcal{Q(R)}$. On the other hand, if $\theta$ is unitary, we can define an
operation, see subsection \ref{_UnitaryGrActsOnStSp}.
\end{remark}

\ 

Let $\mathfrak{B}\in\mathcal{Q}^{ab}\mathcal{(R)}$ be an abelian quasipoint,
and let $E\in\mathfrak{B}$ be an abelian projection. Each $F\in\mathfrak{B}%
_{E}$ is a subprojection of the abelian projection $E$ and hence of the form
$F=QE$, where $Q\in\mathcal{R}$ is a central projection. Then $Q\in
\mathfrak{B}$ holds, so $Q\in\mathfrak{B}\cap\mathcal{C}$. If, conversely,
$Q\in\mathfrak{B}\cap\mathcal{C}$ holds, then $QE\in\mathfrak{B}_{E}$,
therefore we have%
\[
\mathfrak{B}_{E}=\{QE\ |\ Q\in\mathfrak{B}\cap\mathcal{C}\}.
\]
The mapping%
\begin{align*}
\zeta_{E}:\mathcal{R}^{\prime}E  &  \longrightarrow\mathcal{R}^{\prime}%
C_{E},\\
T^{\prime}E  &  \longmapsto T^{\prime}C_{E}%
\end{align*}
is a $\ast$-isomorphism (see Prop. 5.5.5 in \cite{KadRinI97}). Since%
\begin{align*}
PE\wedge QE  &  =PQE=(P\wedge Q)E,\\
PE\vee QE  &  =(P+Q)E-PQE=(P\vee Q)E,
\end{align*}
$\zeta_{E}|_{\mathfrak{B}_{E}}$ is a lattice isomorphism from $\mathfrak{B}%
_{E}$ onto $(\mathfrak{B}\cap\mathcal{C})C_{E}$. $\mathfrak{B}\cap\mathcal{C}$
is a quasipoint of $\mathcal{P(C)}$: obviously, $\mathfrak{B}\cap\mathcal{C}$
is a filter basis $\mathcal{P(C)}$. Let $\beta$ be a quasipoint in
$\mathcal{P(C)}$ containing $\mathfrak{B}\cap\mathcal{C}$. Assume that
$P\in\beta\backslash(\mathfrak{B}\cap\mathcal{C})$. Then $PC_{E}\in
\beta\backslash(\mathfrak{B}\cap\mathcal{C})$ holds, therefore $PE\notin
\mathfrak{B}_{E}$. Then there is some $QE\in\mathfrak{B}_{E}$ such that
$PQE=0$, so $PC_{E}QC_{E}=PQC_{E}=0$, but that contradicts $PC_{E},QC_{E}%
\in\beta$. Thus $\mathfrak{B}\cap\mathcal{C}$ is a quasipoint in
$\mathcal{P(C)}$, and hence%
\[
\zeta_{E}(\mathfrak{B}_{E})=(\mathfrak{B}\cap\mathcal{C})_{C_{E}}.
\]
In this manner, each abelian quasipoint $\mathfrak{B}\in\mathcal{Q}%
^{ab}\mathcal{(R)}$ is assigned a quasipoint $\beta(\mathfrak{B}%
):=\mathfrak{B}\cap\mathcal{C}\in\mathcal{Q(C)}$ of the center $\mathcal{C}$
of $\mathcal{R}$. Moreover, the mapping%
\begin{align*}
\zeta:\mathcal{Q(R)}  &  \longrightarrow\mathcal{Q(C)}\\
\mathfrak{B}  &  \longmapsto\mathfrak{B}\cap\mathcal{C}%
\end{align*}
is surjective, since each quasipoint $\beta\in\mathcal{Q(C)}$ (being a filter
base in $\mathcal{P(R)}$) is contained in some quasipoint $\mathfrak{B}%
\in\mathcal{Q(R)}$.

\ 

Let $\mathfrak{B,}\widetilde{\mathfrak{B}}\in\mathcal{Q}^{ab}\mathcal{(R)}$ be
abelian quasipoints such that%
\[
\beta:=\mathfrak{B}\cap\mathcal{C}=\widetilde{\mathfrak{B}}\cap\mathcal{C}.
\]
Let $E\in\mathfrak{B},\widetilde{E}\in\widetilde{\mathfrak{B}}$ be abelian
projections. Since $C_{E},C_{\widetilde{E}}\in\beta$, $C_{E}C_{\widetilde{E}%
}\in\beta$ holds and $C_{E}C_{\widetilde{E}}E\in\mathfrak{B},C_{E}%
C_{\widetilde{E}}\widetilde{E}\in\widetilde{\mathfrak{B}}$ are abelian
projections with the same central carrier $C_{E}C_{\widetilde{E}}$. Hence,
without loss of generality, one can assume $C_{E}=C_{\widetilde{E}}$. It
follows that $E$ and $\widetilde{E}$ are equivalent (see Prop. 6.4.6 in
\cite{KadRinII97}). Let $\theta\in\mathcal{R}$ be a partial isometry such that
$\theta^{\ast}\theta=E,\theta\theta^{\ast}=\widetilde{E}$, therefore $\theta
E\theta^{\ast}=\widetilde{E}$. It follows that%
\begin{align*}
\theta\mathfrak{B}_{E}\theta^{\ast}  &  =\{\theta QE\theta^{\ast}%
\ |\ Q\in\beta\}=\{Q\theta E\theta^{\ast}\ |\ Q\in\beta\}\\
&  =\{Q\widetilde{E}\ |\ Q\in\beta\}=\mathfrak{B}_{\widetilde{E}},
\end{align*}
so%
\[
\theta_{\mathcal{Q}}(\mathfrak{B})=\widetilde{\mathfrak{B}}.
\]
Conversely, let $\mathfrak{B,}\widetilde{\mathfrak{B}}$ be abelian
quasipoints, and let $\theta\in\mathcal{R}$ be a partial isometry such that
$E:=\theta^{\ast}\theta\in\mathfrak{B},\widetilde{E}:=\theta\theta^{\ast}%
\in\widetilde{\mathfrak{B}}$. From this, $\theta_{\mathcal{Q}}(\mathfrak{B}%
)=\widetilde{\mathfrak{B}}$ as shown. Let $F\in\mathfrak{B}$ be abelian,
$F\leq E$. Then $\theta F$ is a partial isometry from $(\theta F)^{\ast}\theta
F=FEF=F$ to $\theta F\theta^{\ast}\in\widetilde{\mathfrak{B}}$. Since $\theta
F\theta^{\ast}$ is abelian, too, we can assume without loss of generality that
$E$ and $\widetilde{E}$ are abelian. From the definition of $\theta
_{\mathcal{Q}}$, it follows that%
\[
\widetilde{\mathfrak{B}}_{\widetilde{E}}=\theta\mathfrak{B}_{E}\theta^{\ast
}=\{\theta QE\theta^{\ast}\ |\ Q\in\mathfrak{B}\cap\mathcal{C}\}=\{Q\widetilde
{E}\ |\ Q\in\mathfrak{B}\cap\mathcal{C}\}
\]
holds, so%
\[
\{P\widetilde{E}\ |\ P\in\widetilde{\mathfrak{B}}\cap\mathcal{C}%
\}=\{Q\widetilde{E}\ |\ Q\in\mathfrak{B}\cap\mathcal{C}\},
\]
and hence, since $C_{E}=C_{\widetilde{E}}$,%
\begin{align*}
&  \{PC_{E}\ |\ P\in\widetilde{\mathfrak{B}}\cap\mathcal{C}\}=\{QC_{E}%
\ |\ Q\in\mathfrak{B}\cap\mathcal{C}\}\\
\Longleftrightarrow &  (\widetilde{\mathfrak{B}}\cap\mathcal{C})_{C_{E}%
}=(\mathfrak{B}\cap\mathcal{C})_{C_{E}},
\end{align*}
that is, $\widetilde{\mathfrak{B}}\cap\mathcal{C}=\mathfrak{B}\cap\mathcal{C}%
$. Summing up, it is proven that:

\begin{theorem}
\label{TAbQPcQP}Let $\mathcal{R}$ be a von Neumann algebra with center
$\mathcal{C}$. Then the mapping%
\begin{align*}
\zeta:\mathcal{Q(R)}  &  \longrightarrow\mathcal{Q(C)}\\
\mathfrak{B}  &  \longmapsto\mathfrak{B}\cap\mathcal{C}%
\end{align*}
is surjective. If $\mathfrak{B},\widetilde{\mathfrak{B}}\in\mathcal{Q}%
^{ab}(\mathcal{R})$ are two abelian quasipoints, then $\zeta(\mathfrak{B}%
)=\zeta(\widetilde{\mathfrak{B}})$ holds if and only if there is a partial
isometry $\theta\in\mathcal{R}$ with $\theta_{\mathcal{Q}}(\mathfrak{B}%
)=\widetilde{\mathfrak{B}}$.
\end{theorem}

\section{The action of the unitary group on the Stone spectrum $\mathcal{Q(R)}%
$\label{_UnitaryGrActsOnStSp}}

A unitary operator transforms a quasipoint in the obvious way:

\begin{definition}
\label{DUActionOnQP}Let $T\in\mathcal{U(H)}$ be a unitary operator. $T$ acts
on $\mathfrak{B}\in\mathcal{Q(R)}$ $(\mathfrak{R}\subseteq\mathcal{L(H)})$ by%
\[
T.\mathfrak{B}:=\{TET^{\ast}\ |\ E\in\mathfrak{B}\}.
\]

\end{definition}

\begin{lemma}
$T.\mathfrak{B}$ is a quasipoint of the von Neumann algebra $T\mathcal{R}%
T^{\ast}$. If $T\in\mathcal{U(R)}$, then $T.\mathfrak{B}\in\mathcal{Q(R)}$.
\end{lemma}

\begin{proof}
We have%
\[
T(E\wedge F)T^{\ast}\leq TET^{\ast}\wedge TFT^{\ast}.
\]
Moreover,
\[
T^{\ast}(TET^{\ast}\wedge TFT^{\ast})T\leq E\wedge F,
\]
so $T(E\wedge F)T^{\ast}=TET^{\ast}\wedge TFT^{\ast}$. Thus $T.\mathfrak{B}$
is a filter basis and hence contained in some quasipoint $\mathfrak{B}%
^{\prime}\in T\mathcal{R}T^{\ast}$. $T^{\ast}.\mathfrak{B}^{\prime}$ also is a
filterbasis. We have%
\[
T^{\ast}.(T.\mathfrak{B})=\mathfrak{B}\subseteq T^{\ast}.\mathfrak{B}^{\prime
}.
\]
From the maximality of $\mathfrak{B}$, equality holds.
\end{proof}

\section{The Stone spectrum of a type $I_{n}$ von Neumann
algebra\label{_StSpOfTypeI_nAlgebra}}

Let $\mathcal{R}$ be a von Neumann algebra of type $I_{n}$, $n$ finite. We
will show that every quasipoint $\mathfrak{B}\in\mathcal{Q(R)}$ of
$\mathcal{R}$ is abelian, i.e. contains an abelian projection. In order to do
so, we will use the fact that $\mathcal{R}$ is (isomorphic to) a $n\times
n$-\emph{matrix} algebra, albeit with entries from another von Neumann
algebra, the center of $\mathcal{R}$. We regard $\mathcal{R}$ as acting on the
Hilbert module $\mathcal{A}^{n}$, which generalizes the vector space
$\mathbb{C}^{n}$. The abelian projections will be those projecting onto
\textquotedblleft lines\textquotedblright\ of the form $a\mathcal{A}$, where
$\mathcal{A}:=\mathcal{C(R)}$ is the center of $\mathcal{R}$. Of course,
$\mathcal{A}$ is not a field and $\mathcal{A}^{n}$ is not a vector space, so
we cannot use arguments for subspace lattices of finite-dimensional vector
spaces directly (in which case every quasipoint is abelian). But we will
introduce equivalence relations on $\mathcal{A}$ and $\mathcal{A}^{n}$ that
turn them into a field and an $n$-dimensional vector space, respectively, and
show that after taking equivalence classes, enough of the structure remains
intact to allow the conclusion that every quasipoint of $\mathcal{R}$ is
abelian. The intuition from linear algebra carries through. From Thm.
\ref{TAbQPcQP}, we know that the abelian quasipoints can be mapped to
quasipoints of the center of $\mathcal{R}$ via%

\begin{align*}
\xi:\mathcal{Q}^{ab}(\mathcal{R)}  &  \longrightarrow\mathcal{Q(C)},\\
\mathfrak{B}  &  \longmapsto\mathfrak{B}\cap\mathcal{C},
\end{align*}
where two quasipoints $\mathfrak{B},\widetilde{\mathfrak{B}}$ are mapped to
the same quasipoint of the center if and only if there is a partial isometry
$\theta\in\mathcal{R}$ such that $\theta_{\mathcal{Q}}(\mathfrak{B}%
)=\widetilde{\mathfrak{B}}$. Using the fact that $\mathcal{R}$ is a finite
algebra, we can replace partial isometries with unitary operators. This will
allow us to specify the orbits of the unitary group $\mathcal{U(R)}$ acting on
$\mathcal{Q(R)}$ (Thm. \ref{TQPsOfCenterParamUnitaryOrbits}).

\subsection{Hilbert modules and the projections $E_{a}$%
\label{_HilbModsAndProjsE_a}}

It is well known that each type $I_{n}$ von Neumann algebra $\mathcal{R}$ is
$\ast$-isomorphic to $\mathbb{M}_{n}(\mathcal{A})$, the matrix algebra with
entries from $\mathcal{A}=\mathcal{C(R)}$, the center of $\mathcal{R}$ (see
Thm. 6.6.5 in \cite{KadRinII97}). Let $\mathcal{A}^{n}$ be the free right
module over $\mathcal{A}$ consisting of $n$ copies of $\mathcal{A}$. Another
common notation for $\mathbb{M}_{n}(\mathcal{A})$ is $End_{\mathcal{A}%
}(\mathcal{A}^{n})$, the algebra of $\mathcal{A}$-linear endomorphisms of
$\mathcal{A}^{n}$. $\mathbb{M}_{n}(\mathcal{A})$ acts on the Hilbert space
$\widetilde{\mathcal{H}}:=\bigoplus^{n}\mathcal{H}_{\mathcal{A}}$, the
$n$-fold direct sum of $\mathcal{H}_{\mathcal{A}}$, which is the Hilbert space
$\mathcal{A}$ acts on. We will not make use of $\widetilde{\mathcal{H}}$ and
the representation of $\mathbb{M}_{n}(\mathcal{A})$ on it, because we will
regard $\mathbb{M}_{n}(\mathcal{A})$ as an algebra that acts on the
$\mathcal{A}$-module $\mathcal{A}^{n}$ from the left. Elements $a=(a_{1}%
\oplus...\oplus a_{n})^{t}$ of $\mathcal{A}^{n}$ are regarded as column
vectors. The operation of $\mathbb{M}_{n}(\mathcal{A})$ on $\mathcal{A}^{n}$
is a \textquotedblleft matrix$\times$vector\textquotedblright\ operation.
(Since $\mathcal{A}$ is commutative, $\mathcal{A}^{n}$ can be regarded as a
left module as well. The chosen convention fits the natural structure of
$\mathcal{A}^{n}$ as an $\mathbb{M}_{n}(\mathcal{A})$-$\mathcal{A}$-bimodule.)

\ 

$\mathcal{A}^{n}$ has a canonical basis with basis elements%
\[
e_{j}:=(0\oplus...\oplus0\oplus\overset{\underset{\downarrow}{j}}{1}%
\oplus0\oplus...\oplus0)^{t},
\]
where $1$ is the unit of $\mathcal{A}$. With respect to this basis,
$a\in\mathcal{A}^{n}$ is denoted as $a=(a_{1}\oplus...\oplus a_{n})^{t}$. The
sign of transposition will be omitted from now on.\newline

There is an $\mathcal{A}$-valued product defined on $\mathcal{A}^{n}$ such
that $\mathcal{A}^{n}$ becomes a \textbf{Hilbert-}$\mathcal{A}$\textbf{-module%
\index{Hilbert module}}. Since $\mathcal{A}^{n}$ is a right module, the inner
product is $\mathcal{A}$-linear with respect to the \emph{second} variable:%
\begin{align*}
(a|b)  &  =(a_{1}\oplus...\oplus a_{n}|b_{1}\oplus...\oplus b_{n}):=\sum
_{k=1}^{n}a_{k}^{\ast}b_{k},\\
(a|b\alpha)  &  =(a|b)\alpha=\alpha(a|b)=(a\alpha^{\ast}|b)
\end{align*}
for $a,b\in\mathcal{A}^{n},\alpha\in\mathcal{A}$. In the second line the
commutativity of $\mathcal{A}$ was employed. The inner product induces a norm
on $\mathcal{A}^{n}$ by%
\[
|a|:=|(a|a)|^{\frac{1}{2}},
\]
where the norm on the right hand side is the norm on $\mathcal{A}$.

\ 

Let $\Omega:=\mathcal{Q(A)}$ be the Stone spectrum of $\mathcal{A}$. Without
loss of generality, we can assume $\mathcal{A}=C(\Omega)$. Let $(\Omega
_{1},...,\Omega_{n})$ be a partition of $\Omega$ into closed-open sets
$\Omega_{k}\neq\varnothing$, and let $a_{k}:=\chi_{\Omega_{k}}$. Then for all
$\beta\in\Omega$,%
\[
\sum_{k}|a_{k}(\beta)|^{2}=1
\]
and hence $|a_{1}\oplus...\oplus a_{n}|=1$. Moreover, for $k\neq j$,%
\[
(a_{k}|a_{j})=a_{k}^{\ast}a_{j}=0,
\]
that is, the $a_{k}$ are pairwise orthogonal. But (the analogue of)
Pythagoras' theorem does not hold, since for our example one obtains%
\[
|a_{1}\oplus...\oplus a_{n}|^{2}=1<n=\sum_{k=1}^{n}|a_{k}|^{2}.
\]

In general, operators on Hilbert modules are $-$different from those on
Hilbert spaces$-$ not (all) adjointable, which is due to the lack of
self-duality of Hilbert modules, see for example \cite[p 240]{WeO93}. A
mapping $T:E\rightarrow E$ from a Hilbert module $E$ to itself is called
\textbf{adjointable%
\index{Hilbert module!adjointable operator on}} if there is a mapping
$T^{\ast}$ such that%
\[
(Ta|b)=(a|T^{\ast}b)
\]
for all $a,b\in E$. The mapping $T^{\ast}$ is called the \textbf{adjoint of
}$T$. One can show that if $T$ is adjointable, then $T^{\ast}$ is unique,
$T^{\ast\ast}=T$ and both $T$ and $T^{\ast}$ are module maps which are bounded
with respect to the operator norm (\cite[Lemma 15.2.3]{WeO93}).

\ 

For our purpose, we have to characterize the projections in $\mathcal{R}%
\simeq\mathbb{M}_{n}(\mathcal{A)}$.

\begin{lemma}
\label{LProjnInM_n(A)}The elements of $\mathbb{M}_{n}(\mathcal{A})$ are
adjointable, $T=T_{jk}\in\mathbb{M}_{n}(\mathcal{A})$ has adjoint $(T^{\ast
})_{jk}=T_{kj}^{\ast}$.
\end{lemma}

\begin{proof}
Let $a,b\in\mathcal{A}^{n},T=T_{jk}\in\mathbb{M}_{n}(\mathcal{A})$. We have%
\begin{align*}
(Ta|b)  &  =\sum_{k,j}(T_{kj}a_{j})^{\ast}b_{k}=\sum_{k,j}a_{j}^{\ast}%
T_{kj}^{\ast}b_{k}\\
&  =\sum_{j,k}a_{j}^{\ast}(T_{jk})^{\ast}b_{k}=(a|T^{\ast}b).
\end{align*}

\end{proof}

\ 

The adjoint of $T\in\mathbb{M}_{n}(\mathcal{A})$ in the Hilbert module sense
thus is the usual, Hilbert space adjoint of $T$. It follows that the
projections of the von Neumann algebra $\mathbb{M}_{n}(\mathcal{A})$ are the
projections of the algebra $\mathbb{B}(\mathcal{A}^{n})$ of adjointable
operators of the Hilbert module $\mathcal{A}^{n}$. (A projection $P$ is
$\mathcal{A}$-linear, so it is contained in $\mathbb{M}_{n}(\mathcal{A}%
)$.)\newline

We will now introduce projections $E_{a}$ that map from $\mathcal{A}^{n}$ onto
\textquotedblleft lines\textquotedblright\ of the form $a\mathcal{A}$: let
$a\in\mathcal{A}^{n}$ be such that $p:=(a|a)\in\mathcal{A}$ is a projection.
Then%
\[
\forall k\leq n:pa_{k}=a_{k},
\]
since for $\beta\in\mathcal{Q(A)}$ such that $a_{k}(\beta)\neq0$, one obtains
$p(\beta)=%
{\textstyle\sum\nolimits_{j}}
a_{j}^{\ast}a_{j}\neq0$ and hence $p(\beta)=1$, since $p$ is a projection.
Here and in the following, the components $a_{k}\in\mathcal{A}$ of $a$ are
identified with their Gelfand transforms. We get%
\[
\forall\beta\in\mathcal{Q(A)\;\;}\forall k\leq n:(a_{k}p)(\beta)=a_{k}%
(\beta)p(\beta)=a_{k}(\beta),
\]
so $ap=a$.\newline

Now define%
\begin{align*}
E_{a}:\mathcal{A}^{n}  &  \longrightarrow\mathcal{A}^{n},\\
b  &  \longmapsto a(a|b).
\end{align*}
For all $b,c\in\mathcal{A}^{n}$,%
\begin{align*}
(E_{a}b|c)  &  =(a(a|b)|c)=(a|b)^{\ast}(a|c)\\
&  =(b|a)(a|c)=(b|a(a|c))\\
&  =(b|E_{a}c),
\end{align*}
so we have $E_{a}^{\ast}=E_{a}$. Moreover, for all $b\in\mathcal{A}^{n}$,%
\begin{align*}
E_{a}^{2}b  &  =a(a|E_{a}b)=a(a|a(a|b))\\
&  =a(a|a)(a|b)=ap(a|b)\\
&  =a(a|b)=E_{a}b,
\end{align*}
so $E_{a}^{2}=E_{a}$, that is, $E_{a}$ is a projection with
$\operatorname*{im}E_{a}\subseteq a\mathcal{A}$. In fact, equality holds: let
$b=a\alpha\in a\mathcal{A}$. Then%
\[
a(a|b)=a(a|a)\alpha=a\alpha,
\]
therefore $a\mathcal{A}\subseteq\operatorname*{im}E_{a}$. The central carrier
$C_{E_{a}}$ of $E_{a}$ is $I_{n}(a|a)$: obviously, $I_{n}(a|a)$ is a central
projection, since $I_{n}\mathcal{A}\simeq\mathcal{A}$ is the center of
$\mathbb{M}_{n}(\mathcal{A)}$. It holds that%
\begin{align*}
E_{a}b(a|a)  &  =a(a|b)(a|a)=a(a|a)(a|b)\\
&  =a(a|b)=E_{a}b,
\end{align*}
so $C_{E_{a}}\leq I_{n}(a|a)$. Conversely, let $I_{n}q$ be a central
projection such that $qE_{a}=E_{a}q=E_{a}$. Then for all $b\in\mathcal{A}^{n}%
$,%
\begin{align*}
E_{a}b  &  =a(a|b)=qE_{a}b=qa(a|b)\\
&  =a(a|qb)=E_{a}(qb).
\end{align*}
In particular, one obtains
\begin{align*}
&  E_{a}(qa)=E_{a}a\\
\Longleftrightarrow &  a(a|qa)=a(a|a)=ap=a\\
\Longleftrightarrow &  a(a|a)q=a\\
\Longleftrightarrow &  aq=qa=a\\
\implies &  q(a|a)=qp=(a|a)=p\\
\implies &  p\leq q,
\end{align*}
therefore, the central carrier of $E_{a}$ is $C_{E_{a}}=I_{n}p=I_{n}%
(a|a)$.\newline

The $E_{a}$ are of interest, because they are \emph{abelian}
projections:\newline

\begin{lemma}
\label{LE_aAbProj}$E_{a}$ is an abelian projection from $\mathcal{A}^{n}$ onto
$a\mathcal{A}$ with central carrier $I_{n}(a|a)$.
\end{lemma}

\begin{proof}
It only remains to show that $E_{a}$ is abelian. Let $A,B\in\mathbb{M}%
_{n}(\mathcal{A)}$. Then it holds for all $b\in\mathcal{A}^{n}$ that%
\begin{align*}
E_{a}AE_{a}BE_{a}b  &  =E_{a}AE_{a}(Ba)(a|b)=E_{a}(Aa)(a|Ba)(a|b)\\
&  =a(a|Aa)(a|Ba)(a|b)=a(a|Ba)(a|Aa)(a|b)\\
&  =E_{a}(Ba)(a|Aa)(a|b)=E_{a}BE_{a}(Aa)(a|b)\\
&  =E_{a}BE_{a}AE_{a}b,
\end{align*}
so $E_{a}AE_{a}BE_{a}=E_{a}BE_{a}AE_{a}$.
\end{proof}

\ 

$E_{a}$ is a projection in $\mathbb{M}_{n}(\mathcal{A)}$ if $(a|a)$ is a
projection in $\mathcal{A}$. The converse is also true:

\begin{remark}
Let $a\in\mathcal{A}^{n}$ be such that $E_{a}$ is a projection. Then
$(a|a)\in\mathcal{A}$ is a projection.
\end{remark}

\begin{proof}
From $E_{a}^{2}=E_{a}$, $a(a|b)=a(a|E_{a}b)=a(a|a)(a|b)$ for all
$b\in\mathcal{A}^{n}$. For $b=a$,%
\[
a(a|a)=a(a|a)^{2}.
\]
This means that $(a|a)\in\{0,1\}$ holds on the support
\[
supp\ a:=\bigcup_{k\leq n}supp\ a_{k}=supp(a|a).
\]
If $\beta\in\Omega=\mathcal{Q(A)}$ is such that $(a|a)(\beta)\neq0$, then
$a_{k}(\beta)\neq0$ holds for at least one $k\leq n$ and thus $(a|a)(\beta
)=1$. So $(a|a)=1$ holds on $supp(a|a)$ and $(a|a)$ is a projection.
\end{proof}

\ 

If $a_{1},...,a_{n}\in\mathcal{A}$ are projections and $a:=%
{\textstyle\sum\nolimits_{k}}
a_{k}e_{k}$, then $E_{a}$ is a projection if and only if the $a_{k}$ are
pairwise orthogonal, because, according to its definition and the above
remark, $E_{a}$ is a projection if and only if $(a|a)$ is a projection. For
$a:=%
{\textstyle\sum\nolimits_{k}}
a_{k}e_{k}$, we have $(a|a)=%
{\textstyle\sum\nolimits_{k}}
a_{k}$, and this is a projection if and only if the $a_{k}$ are pairwise
orthogonal. Furthermore, $%
{\textstyle\sum\nolimits_{k}}
a_{k}E_{e_{k}}=%
{\textstyle\sum\nolimits_{k}}
E_{a_{k}e_{k}}$ (the $a_{k}$ are projections again): for all $b\in
\mathcal{A}^{n}$, it holds that%
\begin{align*}
(\sum_{k}a_{k}E_{e_{k}})b  &  =\sum_{k}a_{k}e_{k}(e_{k}|b)=\sum_{k}a_{k}%
^{2}e_{k}(e_{k}|b)\\
&  =\sum_{k}a_{k}e_{k}(e_{k}|ba_{k})=\sum_{k}a_{k}e_{k}(a_{k}e_{k}|b)\\
&  =\sum_{k}E_{a_{k}e_{k}}b.
\end{align*}

\begin{lemma}
$A:=%
{\textstyle\sum\nolimits_{k=1}^{n}}
a_{k}E_{e_{k}}$ is a projection if and only if all the $a_{k}$ are
projections. In this case, the central carrier of $A$ is $C_{A}=I_{n}(%
{\textstyle\bigvee\nolimits_{k}}
a_{k})$.
\end{lemma}

\begin{proof}
From $(e_{j}|e_{k})=\delta_{jk}e_{k}$ we get%
\[
\forall c\in\mathcal{A}^{n}:E_{e_{j}}E_{e_{k}}c=e_{j}(e_{j}|E_{e_{k}}%
c)=e_{j}(e_{j}|e_{k})(e_{k}|c)=\delta_{jk}E_{j}c,
\]
so $E_{e_{j}}E_{e_{k}}=\delta_{jk}E_{e_{j}}$. $A:=%
{\textstyle\sum\nolimits_{k=1}^{n}}
a_{k}E_{e_{k}}$ is a projection if and only if $a_{k}^{\ast}=a_{k}$ holds for
all $k\leq n$ and if $A^{2}=A$. Since we have%
\begin{align*}
A^{2}  &  =(\sum_{k=1}^{n}a_{k}E_{e_{k}})^{2}=\sum_{j,k=1}^{n}a_{j}%
a_{k}E_{e_{j}}E_{e_{k}}\\
&  =\sum_{j,k=1}^{n}a_{j}a_{k}\delta_{jk}E_{e_{j}}=\sum_{j=1}^{n}a_{j}%
^{2}E_{e_{j}},
\end{align*}
this holds if and only if all the $a_{j}\in\mathcal{A}$ are projections.
$C_{A}=I_{n}(%
{\textstyle\bigvee\nolimits_{k}}
a_{k})$ holds then, obviously.
\end{proof}

\begin{remark}
\label{RMatrixOfE_a}With respect to the basis $(e_{1},...,e_{n})$ of
$\mathcal{A}^{n}$, $E_{a}$ has the matrix%
\[
(E_{a})_{jk}=(a_{j}a_{k}^{\ast})_{j,k\leq n}.
\]

\end{remark}

\begin{proof}
It holds that%
\begin{align*}
E_{a}e_{k}  &  =a(a|e_{k})=a\sum_{j}a_{j}^{\ast}\delta_{jk}\\
&  =aa_{k}^{\ast}=(\sum_{j}a_{j}e_{j})a_{k}^{\ast}\\
&  =(\sum_{j}a_{j}a_{k}^{\ast})e_{j}.
\end{align*}

\end{proof}

\ 

The projections $E_{a}$ are special cases of the so-called
\textbf{ket-bra-operators%
\index{ket-bra-operators}} (see e.g. \cite[p 71]{GVF01}). These (and their
symbolic Dirac notation) are defined as%
\begin{align*}
|r\rangle\langle s|:E  &  \longrightarrow F,\\
b  &  \longmapsto r(s|b),
\end{align*}
where $E$ and $F$ are Hilbert modules over a $C^{\ast}$-algebra $\mathcal{A}$,
$r\in F$ and $s\in E$. For $E=F$, we have $E_{a}=|a\rangle\langle a|$. There
are some relations among the ket-bra-operators:%
\begin{align*}
r(s|ba)  &  =r(s|b)a,\\
|r\rangle\langle s|^{\ast}  &  =|s\rangle\langle r|,\\
|r\rangle\langle s|\circ|t\rangle\langle u|  &  =|r(t|s)\rangle\langle
u|=|r\rangle\langle u(s|t)|
\end{align*}
hence the finite sums of ket-bra-operators from $E$ to $E$ form a self-adjoint
algebra $End_{\mathcal{A}}^{00}(E)$ contained in $End_{\mathcal{A}}(E)$. The
operators in $End_{\mathcal{A}}^{00}(E)$ are called operators \textbf{of
}$\mathcal{A}$\textbf{-finite rank%
\index{Hilbert module!operator of finite rank}}. The norm closure
$End_{\mathcal{A}}^{0}(E)$ of $End_{\mathcal{A}}^{00}(E)$ contains the
so-called $\mathcal{A}$\textbf{-compact operators%
\index{Hilbert module!a@$\mathcal{A}$-compact operator}}. Clearly,
$End_{\mathcal{A}}^{0}(\mathcal{A}^{n})=\mathbb{M}_{n}(\mathcal{A})$ holds.

\subsection{The modules $a\mathcal{A}$\label{_ModsaA}}

We now turn to the examination of the modules $a\mathcal{A}$ onto which the
$E_{a}$ project. This subsection is quite technical. After showing how to
\textquotedblleft normalize\textquotedblright\ an arbitrary $a\in
\mathcal{A}^{n}\backslash\{0\}$ to $\widetilde{a}$ such that $(\widetilde
{a}|\widetilde{a})$ and $E_{\widetilde{a}}$ are projections with
$a\mathcal{A}=\widetilde{a}\mathcal{A}$, we define the \emph{support} $S(M)$
of a projective submodule $M\subseteq\mathcal{A}^{n}$ and show that there is
an $a\in M$ such that $S(a)=S(M)$ if $M$ is finitely generated. This is used
in the proof that the finitely generated projective submodules of
$\mathcal{A}^{n}$ for which $End_{\mathcal{A}}(M)$ is \emph{abelian} are
exactly those of the form $a\mathcal{A}=\widetilde{a}\mathcal{A}$. To show
this, we also need the fact that $End_{\mathcal{A}}^{0}(P\mathcal{A}%
^{n})=P\mathbb{M}_{n}(\mathcal{A})P$ holds (Lemma 2.18 in \cite{GVF01}).

\ 

All this is a preparation for the following subsection, where quasipoints of
$\mathcal{R}\simeq\mathbb{M}_{n}(\mathcal{A})$ are regarded as families of
projective submodules of $\mathcal{A}^{n}$ rather than families of projections.

\begin{lemma}
Let $a\in\mathcal{A}^{n}$. $a\mathcal{A}$ is a closed submodule of
$\mathcal{A}^{n}$.
\end{lemma}

\begin{proof}
Let $(a\alpha_{n})_{n\in\mathbb{N}}$ be a sequence in $a\mathcal{A}$
converging to $b\in\mathcal{A}^{n}$. As before, let $supp\,a:=%
{\textstyle\bigcup\nolimits_{k\leq n}}
supp\,a_{k}$, then $supp\,a=supp(a|a)$. Without loss of generality, one can
assume that $\alpha_{n}(\beta)=0$ holds for $\beta\neq supp\,a$, since the
sequence $(a\alpha_{n})$ remains unchanged by that. From%
\[
|a\alpha_{n}-a\alpha_{m}|^{2}=|\alpha_{n}-\alpha_{m}|^{2}(a|a)
\]
for all $n,m\in\mathbb{N}$ it follows that $(\alpha_{n})_{n\in\mathbb{N}}$ is
a Cauchy sequence in $\mathcal{A}$. Thus $(\alpha_{n})_{n\in\mathbb{N}}$
converges to some $\alpha\in\mathcal{A}$ and we get%
\[
|a\alpha_{n}-a\alpha|^{2}=|\alpha_{n}-\alpha|^{2}(a|a)\longrightarrow
0\text{\quad for }n\longrightarrow\infty,
\]
that is, $b=a\alpha\in a\mathcal{A}$.
\end{proof}

\ 

Let $M\subseteq\mathcal{A}^{n}$ be some submodule. The \textbf{orthogonal
complement%
\index{Hilbert module!orthogonal complement}} $M^{\bot}$ of $M$ is given by
(see \cite[p 248]{WeO93})%
\[
M^{\bot}:=\{b\in\mathcal{A}^{n}\ |\ \forall a\in M:(b|a)=0\}.
\]
For $M=a\mathcal{A}$, we obtain%
\begin{align*}
(a\mathcal{A)}^{\bot}  &  =\{b\in\mathcal{A}^{n}\ |\ \forall\alpha
\in\mathcal{A}:(b|a\alpha)=0\}\\
&  =\{b\in\mathcal{A}^{n}\ |\ \forall\alpha\in\mathcal{A}:(b|a)\alpha=0\}\\
&  =\{b\in\mathcal{A}^{n}\ |\ (b|a)=0\}.
\end{align*}
Obviously, $M\cap M^{\bot}=0$ holds for any submodule $M$. A submodule $M$ is
called \textbf{complementable} if $M\oplus M^{\bot}=\mathcal{A}^{n}$ holds.
One can show that $M$ is complementable if and only if it is the image of some
projection (Cor. 15.3.9 in \cite{WeO93}).\newline

Subsequently, it will be demonstrated how to normalize an arbitrary
$a\in\mathcal{A}^{n}\backslash\{0\}$ to $\widetilde{a}$ such that
$(\widetilde{a}|\widetilde{a})$ is a projection in $\mathcal{A}$ and
$E_{\widetilde{a}}$ is the projection in $\mathbb{M}_{n}(\mathcal{A)}$ onto
$a\mathcal{A}=\widetilde{a}\mathcal{A}$. Let $a\in\mathcal{A}^{n}%
\backslash\{0\}$ be such that $(a|a)$ is not a projection. For $n\in
\mathbb{N}$, let%
\[
A_{n}:=\{\beta\in\Omega\ |\ (a|a)(\beta)>\frac{1}{n}\}.
\]
$A_{n}$ is open and $\overline{A_{n}}$ is open and closed, since $\Omega$ is
extremely disconnected. For an appropriate $n_{0}\in\mathbb{N}$, $A_{n}%
\neq\varnothing$ holds for all $n\geq n_{0}$.

\ 

Let $\alpha_{n}\in\mathcal{A}$ be defined by%
\[
\alpha_{n}(\beta):=\left\{
\begin{tabular}
[c]{ll}%
$(a|a)(\beta)^{-\frac{1}{2}}$ & for $\beta\in\overline{A_{n}}$\\
$0$ & for $\beta\notin\overline{A_{n}}$%
\end{tabular}
\ \ \ \ \ \right.  .
\]
Then $a\alpha_{n}\in a\mathcal{A}$ and%
\[
(a\alpha_{n}|a\alpha_{n})=(a|a)\alpha_{n}^{2}=\left\{
\begin{tabular}
[c]{ll}%
$1$ & on $\overline{A_{n}}$\\
$0$ & on $\Omega\backslash\overline{A_{n}}$%
\end{tabular}
\ \ \ \ \ \right.  ,
\]
therefore $(a\alpha_{n}|a\alpha_{n})$ is a projection in $\mathcal{A}$. Let
$E_{n}:=E_{a\alpha_{n}}$ be the projection onto $a\alpha_{n}\mathcal{A}$ given
by $a\alpha_{n}$. According to the definition of $\alpha_{n}$, $\alpha
_{n}\mathcal{A}$ is the closed ideal%
\[
\alpha_{n}\mathcal{A}=\{\alpha\in\mathcal{A}\ |\ supp\,\alpha\subseteq
\overline{A_{n}}\}=\chi_{\overline{A_{n}}}\mathcal{A}.
\]
For these ideals it holds that%
\[
\forall n\in\mathbb{N}:\alpha_{n}\mathcal{A}\subseteq\alpha_{n+1}\mathcal{A}.
\]
Moreover, with $S(a):=supp\,a\subseteq\Omega=\mathcal{Q(A)}$, it holds that
$\overline{\bigcup_{n}\alpha_{n}\mathcal{A}}=\mathcal{A}\chi_{S(a)}$: the
inclusion \textquotedblleft$\subseteq$\textquotedblright\ is clear from
$\overline{A_{n}}\subseteq S(a)$. Let $b\in\mathcal{A}^{n}$ be such that
$S(b)\subseteq S(a)$, with no loss of generality $b\geq0$. Then $(b\chi
_{\overline{A_{n}}})_{n\geq n_{0}}$ is increasing monotonously. One gets%
\[
\sup_{n\in\mathbb{N}}b\chi_{\overline{A_{n}}}=b\sup_{n\in\mathbb{N}}%
\chi_{\overline{A_{n}}}.
\]
From%
\[
\bigcup_{n}A_{n}=\{\beta\in\Omega\ |\ (a|a)(\beta)>0\},
\]
we have%
\[
\overline{\bigcup_{n}\overline{A_{n}}}=\overline{\bigcup_{n}A_{n}}=S(a)
\]
and thus $\sup_{n}\chi_{\overline{A_{n}}}=\chi_{S(a)}$, therefore%
\[
b=\sup_{n}b\chi_{\overline{A_{n}}}.
\]
According to Dini's theorem, $b\chi_{\overline{A_{n}}}$ converges uniformly to
$b$, and we have shown that $\overline{\bigcup_{n}\alpha_{n}\mathcal{A}%
}=\mathcal{A}\chi_{S(a)}$ holds.\newline

The central carrier of the projection $E_{n}$ is $\chi_{\overline{A_{n}}}$.
$E_{n}\leq E_{n+1}$ holds, since for all $b\in\mathcal{A}^{n}$,%
\begin{align*}
E_{n}E_{n+1}b  &  =a\alpha_{n}(a\alpha_{n}|E_{a\alpha_{n+1}}b)=a\alpha
_{n}(a\alpha_{n}|a\alpha_{n+1})(a\alpha_{n+1}|b)\\
&  =a\alpha_{n}(a\alpha_{n+1}|a\alpha_{n+1})(a\alpha_{n}|b)=a\alpha_{n}%
\chi_{\overline{A_{n+1}}}(a\alpha_{n}|b)\\
&  =a\alpha_{n}(a\alpha_{n}|b)=E_{n}b,
\end{align*}
where we used $S(a\alpha_{n})=\overline{A_{n}}\subseteq\overline{A_{n+1}}$ in
the penultimate step. Let $E:=\bigvee\nolimits_{n\in\mathbb{N}}E_{n}$. The
image of $E_{n}$ is $a\mathcal{A}\chi_{\overline{A_{n}}}$, so%
\[
\operatorname*{im}E=a\mathcal{A}\chi_{S(a)}=a\mathcal{A}.
\]
When defining $E$, one cannot simply assume the properties of the Hilbert
space situation.\newline

The sesquilinear form%
\begin{align*}
\mathcal{A}^{n}\times\mathcal{A}^{n}  &  \longrightarrow\mathcal{A},\\
(b,c)  &  \longmapsto(E_{n}b|c)
\end{align*}
can be written as%
\begin{align*}
(E_{n}b|c)  &  =(a\alpha_{n}(a\alpha_{n}|b)|c)=\alpha_{n}^{\ast}%
(a|c)(a\alpha_{n}|b)^{\ast}\\
&  =\alpha_{n}^{\ast}(a|c)(b|a)\alpha_{n}=\chi_{\overline{A_{n}}%
}(a(a|a)^{-\frac{1}{2}}|c)(b|a(a|a)^{-\frac{1}{2}}),
\end{align*}
and for $n\rightarrow\infty$, the right hand side converges to%
\[
(Eb|c):=\chi_{S(a)}(a(a|a)^{-\frac{1}{2}}|b)(a(a|a)^{-\frac{1}{2}}|c).
\]
Here we use the fact that every bounded continuous function $f$ on an open
dense subset $G\subseteq\Omega$ of the Stone space $\Omega$ can be extended to
a continuous function $\widetilde{f}$ on the whole of $\Omega$ (see Cor.
III.1.8 in \cite{TakI02}): the support $S(a)=supp\,a\subseteq\Omega$ is open
and closed, $G(a):=\{\beta\in\Omega\ |\ a(\beta)\neq0\}$ is open and dense in
$S(A)$ and the mapping $a(a|a)^{-\frac{1}{2}}:G(a)\rightarrow\mathbb{C}^{n}$
is continuous und bounded, therefore it can be extended to a continuous
mapping $S(a)\rightarrow\mathbb{C}^{n}$, which will also be denoted by
$a(a|a)^{-\frac{1}{2}}$. Define%
\[
\widetilde{a}(\beta):=\left\{
\begin{tabular}
[c]{ll}%
$a(a|a)^{-\frac{1}{2}}$ & $\text{for }\beta\in S(a)$\\
$0$ & $\text{for }\beta\in\Omega\backslash S(a)$%
\end{tabular}
\ \ \ \ \ \ \ \ \right.
\]
and from this, $E_{\widetilde{a}}$. Then it holds for all $b,c\in
\mathcal{A}^{n}$ that%
\[
(E_{\widetilde{a}}b|c)=(\widetilde{a}|c)(\widetilde{a}|b)=\chi_{S(a)}%
(a(a|a)^{-\frac{1}{2}}|c)(a(a|a)^{-\frac{1}{2}}|b).
\]
Obviously, $E_{\widetilde{a}}$ is the same as the limit $E$ of the projections
$E_{n}$ defined above. We showed both ways of the definition, because we will
need the sets $A_{n}$ for the following lemma.

\ 

Let $a\in\mathcal{A}^{n}\backslash\{0\}$, and let $\widetilde{a}\in
\mathcal{A}^{n}$ be as defined above. Then%
\[
(\widetilde{a}|\widetilde{a})=(a|a)^{-1}(a|a)=1
\]
holds on $G(a)$ and hence also on $S(a)$. Thus $(\widetilde{a}|\widetilde{a})$
is a projection in $\mathcal{A}$.

\begin{lemma}
\label{La(wave)A=aA}For the closed submodules, we have $\widetilde
{a}\mathcal{A}=a\mathcal{A}$.
\end{lemma}

\begin{proof}
From $\overline{%
{\textstyle\bigcup\nolimits_{n}}
\mathcal{A\chi}_{\overline{A_{n}}}}=\mathcal{A}\chi_{S(a)}$, one gets
$a\mathcal{A=}\overline{%
{\textstyle\bigcup\nolimits_{n}}
a\mathcal{A}\chi_{\overline{A_{n}}}}$. Moreover, it holds that%
\[
a\mathcal{A}\chi_{\overline{A_{n}}}=a(a|a)^{-\frac{1}{2}}\mathcal{A}%
\chi_{\overline{A_{n}}}=\widetilde{a}\mathcal{A}\chi_{\overline{A_{n}}}.
\]
It follows from this and $S(a)=S(\widetilde{a})$ that%
\[
a\mathcal{A}=\overline{%
{\textstyle\bigcup\nolimits_{n}}
a\mathcal{A}\chi_{\overline{A_{n}}}}=\overline{%
{\textstyle\bigcup\nolimits_{n}}
\widetilde{a}\mathcal{A}\chi_{\overline{A_{n}}}}=\widetilde{a}\mathcal{A}.
\]

\end{proof}

\ 

Of course, $\widetilde{a}=a$ if $(a|a)$ is a projection in $\mathcal{A}$, so
one obtains

\begin{proposition}
\label{PE_a(wave)AbProjOntoaA}For every $a\in\mathcal{A}^{n}$, $E_{\widetilde
{a}}$ is an abelian projection with image $a\mathcal{A}$.
\end{proposition}

$E_{\widetilde{a}}$ is the unique projection from $\mathcal{A}^{n}$ onto
$a\mathcal{A}$: let $Q:\mathcal{A}^{n}\rightarrow a\mathcal{A}$ be a
projection. Then%
\[
\forall c\in(a\mathcal{A})^{\bot}:(Qc|Qc)=(Qc|c)=0,
\]
so $Q|_{(a\mathcal{A})^{\bot}}=0$. Let $Qa=a\alpha$. From $Q^{2}=Q$,
$a\alpha^{2}=a\alpha$, therefore $\alpha(\beta)\in\{0,1\}$ on $S(a)$. Since
$Q\mathcal{A}^{n}=a\mathcal{A}$, it follows that $\alpha=1$ holds on $S(a)$,
so $Qa=a$. Let $b\in\mathcal{A}^{n},b=a\gamma+a^{\prime}$ be such that
$\gamma\in\mathcal{A}$ and $a^{\prime}\in(a\mathcal{A})^{\bot}$. Such a
decomposition exists, since $a\mathcal{A}$ is the image of a projection and
hence a complementable submodule of $\mathcal{A}^{n}$. We have $Qb=a\gamma$
and%
\begin{align*}
E_{\widetilde{a}}b  &  =\widetilde{a}(\widetilde{a}|b)=\widetilde
{a}(\widetilde{a}|a)\gamma=\widetilde{a}(\widetilde{a}|(a|a)^{\frac{1}{2}%
}\widetilde{a})\gamma\\
&  =\widetilde{a}(a|a)^{\frac{1}{2}}\gamma=a\gamma=Qb,
\end{align*}
so $Q=E_{\widetilde{a}}$.\newline

Next, we will define the support of a submodule $M\subseteq\mathcal{A}^{n}$.
For this, we will need a

\begin{remark}
Let $U,V\subseteq\Omega$ be open. If $U\cap V=\varnothing$, then $\overline
{U}\cap\overline{V}=\varnothing$.
\end{remark}

\begin{proof}
Let $\overline{U}\cap\overline{V}\neq\varnothing$ and $\beta\in\overline
{U}\cap\overline{V}$. Since $\overline{U}$ is open, $\overline{U}\cap
V\neq\varnothing$. Since $V$ is open, $U\cap V\neq\varnothing$.
\end{proof}

\ 

Let $M$ be a \emph{projective submodule}%
\index{Hilbert module!projective submodule of} of $\mathcal{A}^{n}$, i.e.
$M=P\mathcal{A}^{n}$ for some projection $P\in\mathbb{M}_{n}(\mathcal{A})$
(see p 89 in \cite{GVF01}). Then $M$ is finitely generated.%
\[
ann(M):=\{\alpha\in\mathcal{A}\ |\ M\alpha=0\},
\]
the \textbf{annihilator%
\index{Hilbert module!annihilator of projective submodule} of }$M$, is a
closed ideal in $\mathcal{A}$. Let $M$ be finitely generated, $\{g_{1}%
,...,g_{r}\}\subseteq M$ a system of generators. For $\alpha$ $\in\mathcal{A}%
$, let $P(\alpha):=\{\beta\in\Omega\ |\ \alpha(\beta)\neq0\}$ and for
$a\in\mathcal{A}^{n}$ let $P(a):=\{\omega\in\Omega\ |\ a(\omega)\neq0\}$. Then%
\[
\alpha a=0\Longleftrightarrow P(\alpha)\cap P(a)=\varnothing,
\]
and thus, according to the above remark,%
\[
\alpha a=0\Longleftrightarrow S(\alpha)\cap S(a)=\varnothing.
\]
Let%
\[
P(M):=\{\beta\in\Omega\ |\ \exists a\in M:\beta\in P(a)\}.
\]
$S(M):=\overline{P(M)}$ is called the \textbf{support%
\index{Hilbert module!support of projective submodule} of }$M$. Since
$\{g_{1},...,g_{r}\}$ is a system of generators of $M$, $P(M)=%
{\textstyle\bigcup\nolimits_{k\leq r}}
P(g_{k})$ holds and hence $S(M)=%
{\textstyle\bigcup\nolimits_{k\leq r}}
S(g_{k})=%
{\textstyle\bigcup\nolimits_{a\in M}}
S(a)$. The set $S(M)$ is open and closed. $\alpha$ $\in ann(M)$ holds if and
only if $\alpha g_{k}=0$ for all $k\leq r$, therefore%
\begin{align*}
&  \alpha\in ann(M)\\
\Longleftrightarrow &  S(\alpha)\cap S(g_{k})=\varnothing\text{\quad}(k\leq
r)\\
\Longleftrightarrow &  S(\alpha)\cap S(M)=\varnothing\\
\Longleftrightarrow &  S(\alpha)\subseteq\Omega\backslash S(M)\\
\Longleftrightarrow &  \alpha|_{S(M)}=0.
\end{align*}
This shows:

\begin{remark}
$ann(M)$ is the vanishing ideal of the closed-open set $S(M)$.
\end{remark}

\begin{lemma}
$I_{n}\chi_{S(M)}$ is the central carrier of the projection $P_{M}%
\in\mathbb{M}_{n}(\mathcal{A)}$ from $\mathcal{A}^{n}$ onto $M$.
\end{lemma}

\begin{proof}
Since $\chi_{S(M)}=1$ on $S(a)$ for all $a\in M$, we have $\chi_{S(M)}a=a$ for
all $a\in M$ and hence $\chi_{S(M)}P_{M}=P_{M}$, i.e. $C_{P(M)}\leq I_{n}%
\chi_{S(M)}$. Let $p\in\mathcal{A}$ be a projection such that $I_{n}p\geq
P_{M}$. Then $pa=pP_{M}a=P_{M}a=a$ for all $a\in M$, so $p=1$ on $S(a)$, that
is, $p\geq\chi_{S(a)}$. It follows that $p\geq\chi_{S(M)}$.
\end{proof}

\begin{lemma}
Let $M$ be a finitely generated, projective submodule of $\mathcal{A}^{n}$,
and let $a,b\in M$. Then there is a $c\in M$ such that $S(a)\cup S(b)\subseteq
S(c)$.
\end{lemma}

\begin{proof}
Regard the decomposition $b=a\alpha+a^{\prime}$ with $\alpha\in\mathcal{A}$
and $a^{\prime}\in(a\mathcal{A)}^{\bot}$. Such a decomposition always exists,
since $a\mathcal{A}$ is a projective and hence complementable submodule of
$\mathcal{A}^{n}$. Therefore each $b\in\mathcal{A}^{n}$ can be decomposed, in
particular each $b\in M\subseteq\mathcal{A}^{n}$. We have $a^{\prime
}=b-a\alpha\in M$, so $a+a^{\prime}\in M$ and%
\begin{align*}
(b|b)  &  =\alpha^{\ast}\alpha(a|a)+(a^{\prime}|a^{\prime}),\\
(a+a^{\prime}|a+a^{\prime})  &  =(a|a)+(a^{\prime}|a^{\prime}),
\end{align*}
therefore
\begin{align*}
(a|a)(\beta)>0  &  \implies(a+a^{\prime}|a+a^{\prime})(\beta)>0,\\
(b|b)(\beta)>0  &  \implies(a|a)(\beta)>0\text{ or }(a^{\prime}|a^{\prime
})(\beta)>0\\
&  \implies(a+a^{\prime}|a+a^{\prime})(\beta)>0.
\end{align*}
It follows that $S(a|a):=supp(a|a)\subseteq S(a+a^{\prime}|a+a^{\prime})$ and
$S(b|b)\subseteq S(a+a^{\prime}|a+a^{\prime})$, that is, $S(a)\cup
S(b)\subseteq S(a+a^{\prime})$.
\end{proof}

\begin{corollary}
\label{CS(M)=S(a)}Let $M$ be as before. Then there is an $a\in M$ such that
$S(a)=S(M)$. $a$ can be chosen such that $(a|a)$ is a projection.
\end{corollary}

\begin{proof}
Let $u_{1},...,u_{r}\in M$ with $M=u_{1}\mathcal{A}+...+u_{r}\mathcal{A}$.
According to the lemma above, there is an\ $a\in M$ such that $S(M)=%
{\textstyle\bigcup\nolimits_{k\leq r}}
S(u_{k})\subseteq S(a)\subseteq S(M)$, so $S(M)=S(a)$. Since $\widetilde
{a}\mathcal{A}=a\mathcal{A}\subseteq M$, $\widetilde{a}\in M$, $(\widetilde
{a}|\widetilde{a})$ is a projection and $S(\widetilde{a})=S(a)$.
\end{proof}

\begin{lemma}
\label{LPropertiesOfKetBraOps}For ket-bra-operators $|b\rangle\langle
a|,|v\rangle\langle u|\in\mathbb{M}_{n}(\mathcal{A})$ it holds that
$|b\rangle\langle a|\circ|v\rangle\langle u|=(a|v)|b\rangle\langle u|$
($a,b,u,v\in M\subseteq\mathcal{A}^{n}$).
\end{lemma}

\begin{proof}
Let $c\in M$. Then%
\begin{align*}
|b\rangle\langle a|\circ|v\rangle\langle u|(c)  &  =b(a||v\rangle\langle
u|(c))=b(a|v(u|c))\\
&  =b(a|v)(u|c)=(a|v)|b\rangle\langle u|(c).
\end{align*}

\end{proof}

\ 

Let $M=P\mathcal{A}^{n}$ be a projective submodule. Lemma 2.18 in \cite{GVF01}
tells us that $End_{\mathcal{A}}^{0}(M)=P\mathbb{M}_{n}(\mathcal{A})P$, hence
if $P$ is abelian, so is $End_{\mathcal{A}}^{0}(P\mathcal{A}^{n})$ and, in
particular, $End_{\mathcal{A}}^{00}(P\mathcal{A}^{n})$ is abelian. We will now
characterize the finitely generated projective submodules $M\subseteq
\mathcal{A}^{n}$ for which $End_{\mathcal{A}}^{00}(M)$, the set of
$\mathcal{A}$-linear mappings of $\mathcal{A}$-finite rank from $M$ to itself,
is abelian. Let $M$ be such a module. Let $u_{1},...,u_{r}\in M$ be generators
of $M$ such that $(u_{k}|u_{k})\in\mathcal{A}$ is a projection for all $k\leq
r$. Moreover, let $a\in M$ be such that $(a|a)$ is a projection and
$S(a)=S(M)$ holds. Then%
\[
|u_{k}\rangle\langle a|\circ|a\rangle\langle u_{k}|=(a|a)|u_{k}\rangle\langle
u_{k}|
\]
is a projection from $M$ \emph{onto} $(a|a)u_{k}\mathcal{A}$, and%
\[
|a\rangle\langle u_{k}|\circ|u_{k}\rangle\langle a|=(u_{k}|u_{k}%
)|a\rangle\langle a|
\]
is a projection from $M$ \emph{onto} $(u_{k}|u_{k})a\mathcal{A}$. Since
$End_{\mathcal{A}}^{00}(M)$ is abelian by assumption, we get%
\[
\forall k\leq r:(a|a)u_{k}\mathcal{A}=(u_{k}|u_{k})a\mathcal{A},
\]
and since $S(a|a)=S(M)$, it holds that%
\[
\forall k\leq r:u_{k}\mathcal{A}=(u_{k}|u_{k})a\mathcal{A}\subseteq
a\mathcal{A}.
\]
Thus it follows that%
\[
M=\sum_{k}u_{k}\mathcal{A}\subseteq a\mathcal{A}\subseteq M,
\]
that is, $M=a\mathcal{A}$, so $M$ is simply generated. Summing up, we have shown:

\begin{proposition}
\label{PProjAbIffP=E_aEtc}A projection $P\in\mathbb{M}_{n}(\mathcal{A}%
)=End_{\mathcal{A}}(\mathcal{A}^{n})$ is abelian if and only if there is an
$a\in\mathcal{A}^{n}$ such that $P=E_{a}$. Then $(a|a)\in\mathcal{A}$ is a
projection and $I_{n}(a|a)$ is the central carrier of $E_{a}$. For
$a\in\mathcal{A}^{n}$, $E_{\widetilde{a}}$ \ is the unique projection onto the
simply generated submodule $a\mathcal{A}$ of $\mathcal{A}^{n}$. This submodule
is projective. $a\mathcal{A}$ is free if and only if $S(a)=\Omega
=\mathcal{Q(A)}$, i.e. if $E_{\widetilde{a}}$ has central carrier $I_{n}$.
\end{proposition}

\subsection{The equivalence relation on $\mathcal{A}^{n}$ and abelian
quasipoints of $\mathcal{P}(\mathbb{M}_{n}(\mathcal{A}))$%
\label{_EquivRelOnAnAndAbQPs}}

In the following, we will regard quasipoints $\mathfrak{B}$ of $\mathcal{P}%
(\mathbb{M}_{n}(\mathcal{A}))$ as families of projective submodules of
$\mathcal{A}^{n}$ rather than families of projections. The filter basis
property of quasipoints will allow us to define a certain notion of germs on
$\mathcal{A}$ and $\mathcal{A}^{n}$ and reduce the situation to that of finite
dimensional vector spaces. The results known from that simple part of the
theory (see Prop. \ref{PQPsOfTypeI_nFactors}) easily show that each quasipoint
$\mathfrak{B}\in\mathcal{Q}(\mathbb{M}_{n}(\mathcal{A}))$ contains a submodule
of the form $a_{0}\mathcal{A}$ and hence is abelian.\newline

So let $\mathfrak{B}\in\mathcal{Q}(\mathbb{M}_{n}(\mathcal{A))}$ a quasipoint,
regarded as a family of projective submodules $M$ of $\mathcal{A}^{n}$. We
have%
\[
P_{M\cap N}=P_{M}\wedge P_{N},
\]
if $M$, $N$ as well as $M\cap N$ are algebraically finitely generated (see
Sec. 15.4 of \cite{WeO93}).\newline

Let $P_{M}\wedge P_{N}$ be the minimum of the projections $P_{M},P_{N}%
\in\mathbb{M}_{n}(\mathcal{A})$. Then%
\[
P_{M}\wedge P_{N}=P_{K},
\]
where $K$ is the largest finitely generated closed submodule of $\mathcal{A}%
^{n}$ such that $K\subseteq M\cap N$ holds. We will use the notation
$K=M\wedge N$.\newline

As mentioned above, $I_{n}\mathcal{A}$ is the center of $\mathbb{M}%
_{n}(\mathcal{A)}$.%
\[
\beta:=\mathfrak{B}\cap I_{n}\mathcal{A}%
\]
is a quasipoint of $\mathcal{A}$, and it holds that%
\[
\beta\simeq\{\chi_{S(M)}\ |\ M\in\mathfrak{B}\}.
\]
This can also be expressed in the following way:%
\begin{align*}
&  \forall M\in\mathfrak{B}:\beta\in\bigcup_{a\in M}S(a)\\
\Longleftrightarrow &  \forall M\in\mathfrak{B}:\chi_{S(M)}(\beta)=1.
\end{align*}
Notice the double role of $\beta$:\ on the one hand as an element of the Stone
spectrum (that is, the Gelfand spectrum) $\mathcal{Q(A)}$ of the center
$I_{n}\mathcal{A}\simeq\mathcal{A}$ of $\mathbb{M}_{n}(\mathcal{A)}$, on the
other hand as a collection of characteristic functions $\chi_{S(M)}$ on this
Stone spectrum.

\ 

We now define an equivalence relation on $\mathcal{A}^{n}$ that amounts to
taking germs:

\begin{definition}
\label{DEquivRelOnAn}Let $n\in\mathbb{N}$. Two elements $a,b$ of
$\mathcal{A}^{n}$ are called \textbf{equivalent at the quasipoint%
\index{Hilbert module!equivalence at quasipoint} }$\beta\in\mathcal{Q(A)}$, if
there is a $p\in\beta$ such that $pa=pb$. Notation: $a\sim_{\beta}b$.
\end{definition}

$\sim_{\beta}$ really is an equivalence relation: symmetry and reflexivity are
obvious. Let $a\sim_{\beta}b,b\sim_{\beta}c$, then there are $p,q\in\beta$
such that $pa=pb,qb=qc\Rightarrow pq\in\beta$ and $pqa=pqb=pqc$, therefore
$a\sim_{\beta}c$. Let $[a]_{\beta}$ be the equivalence class of $a\in
\mathcal{A}^{n}$, and let $[\mathcal{A}^{n}]_{\beta}:=\{[a]_{\beta}%
\ |\ a\in\mathcal{A}^{n}\}$.

\begin{theorem}
\label{T[A]Field[An]VectorSpace}(i) $[\mathcal{A}]_{\beta}$ is a
field,\newline(ii) $[\mathcal{A}^{n}]_{\beta}$ is an $n$-dimensional vector
space over $[\mathcal{A]}_{\beta}.$
\end{theorem}

\begin{proof}
Let $a,b\in\mathcal{A}^{n}$. Then%
\[
\lbrack a]_{\beta}+[b]_{\beta}:=[a+b]_{\beta}%
\]
is well defined, and also, for $\alpha\in\mathcal{A}$,%
\[
\lbrack a]_{\beta}[\alpha]_{\beta}:=[a\alpha]_{\beta}%
\]
is well defined: from $a\sim_{\beta}a^{\prime},b\sim_{\beta}b^{\prime}%
,\alpha\sim_{\beta}\alpha^{\prime}$ it follows that there exist $p,q,r\in
\beta$ such that%
\begin{align*}
&  pa=pa^{\prime},qb=qb^{\prime},r\alpha=r\alpha^{\prime}\\
\Longrightarrow &  pq(a+b)=pqa+pqb=pqa^{\prime}+pqb^{\prime}=pq(a^{\prime
}+b^{\prime})
\end{align*}
and%
\[
pr(a\alpha)=pa(r\alpha)=pa^{\prime}(r\alpha^{\prime})=pr(a^{\prime}%
\alpha^{\prime}).
\]
Furthermore,%
\[
\lbrack\alpha]_{\beta}[\gamma]_{\beta}:=[\alpha\gamma]_{\beta}%
\]
defines a multiplication ($\alpha,\gamma\in\mathcal{A}$): it holds that%
\begin{align*}
&  a\sim_{\beta}0\\
\Longleftrightarrow &  \exists p\in\beta:pa=0\\
\Longleftrightarrow &  \exists p\in\beta:P(p)\cap P(a)=\varnothing\\
\Longleftrightarrow &  \exists p\in\beta:S(p)\cap S(a)=\varnothing\\
\Longleftrightarrow &  \beta\notin S(a).
\end{align*}
Let $\alpha\in\mathcal{A}$ be such that $[\alpha]_{\beta}\neq0$. Then
$\alpha^{\ast}\alpha\geq\varepsilon>0$ holds on a closed-open neighbourhood
$W$ of $\beta$ in $\Omega=\mathcal{Q(A)}$, since from $p\alpha\neq0$ for all
$p\in\beta$ and $p(\beta)=1$, it follows that $\alpha(\omega)\neq0$ holds on a
neighbourhood of $\beta$. Thus $\chi_{W}\alpha$ is invertable on $W$ and there
is a $\gamma\in C(\Omega)$ such that $S(\gamma)=W$ and $\chi_{W}\alpha
\gamma=\chi_{W}$, that is, $\alpha\gamma\sim_{\beta}1$, so $[\alpha]_{\beta
}[\gamma]_{\beta}=1$. Since obviously the algebraic rules for multiplication
and addition are fulfilled, it follows that $[\mathcal{A}]_{\beta}$ is a field
and $[\mathcal{A}^{n}]_{\beta}$ is a vector space over $[\mathcal{A}]_{\beta}%
$. $([e_{1}]_{\beta},...,[e_{n}]_{\beta})$ is a basis of $[\mathcal{A}%
^{n}]_{\beta}$: let $a\in\mathcal{A}^{n},a=%
{\textstyle\sum\nolimits_{k=1}^{n}}
e_{k}a_{k}$, then $[a]_{\beta}=%
{\textstyle\sum\nolimits_{k}}
[e_{k}]_{\beta}[a_{k}]_{\beta}$; if $%
{\textstyle\sum\nolimits_{k}}
[e_{k}]_{\beta}[\gamma_{k}]_{\beta}=0$, then there is some $p\in\beta$ such
that%
\[
0=p\sum_{k}e_{k}\gamma_{k}=\sum_{k}e_{k}(p\gamma_{k}),
\]
so $p\gamma_{k}=0$ for all $k$, that is, $[\gamma_{k}]_{\beta}=0$ for all $k$.
Thus we have%
\[
\dim_{[\mathcal{A}]_{\beta}}[\mathcal{A}^{n}]_{\beta}=n.
\]

\end{proof}

\ 

Let $M\subseteq\mathcal{A}^{n}$ be a submodule. Then $[M]_{\beta}%
:=\{[\alpha]_{\beta}\ |\ a\in M\}$ is a subspace of $[\mathcal{A}^{n}]_{\beta
}$. If $N\subseteq\mathcal{A}^{n}$ is another submodule, then

\begin{lemma}
$[M\cap N]_{\beta}=[M]_{\beta}\cap\lbrack N]_{\beta}$.
\end{lemma}

\begin{proof}
The inclusion \textquotedblleft$\subseteq$\textquotedblright\ is trivial. Let
$[a]_{\beta}\in\lbrack M]_{\beta}\cap\lbrack N]_{\beta}$. Then $[a]_{\beta
}=[b]_{\beta}$ holds for some $b\in N$, so $pa=pb$ for a $p\in\beta$ and hence
$pa\in M\cap N$. Since $p\sim_{\beta}1$ in $\mathcal{A}$ ($p$ is a
projection), it follows that $[a]_{\beta}=[pa]_{\beta}\in\lbrack M\cap
N]_{\beta}$.
\end{proof}

\begin{corollary}
$M,N\in\mathfrak{B}\Rightarrow\lbrack M\wedge N]_{\beta}\subseteq\lbrack
M]_{\beta}\cap\lbrack N]_{\beta}$.
\end{corollary}

\begin{proof}
This follows from the lemma and $M\wedge N\subseteq M\cap N$.
\end{proof}

\ 

Let $M\in\mathfrak{B}$, then $[M]_{\beta}\neq0$: assume that $[M]_{\beta}=0$.
Then%
\begin{align*}
&  \forall a\in M\quad\exists p_{a}\in\beta:p_{a}a=0\\
\implies &  \forall a\in M:\beta\notin S(a),
\end{align*}
but $\beta\in S(M)=\bigcup\nolimits_{a\in M}S(a)$, since $M\in\mathfrak{B}$.
Thus we get

\begin{remark}
\label{RB_bFilterBase}$\mathfrak{B}_{\beta}:=\{[M]_{\beta}\ |\ M\in
\mathfrak{B}\}\subseteq\mathbb{L}([\mathcal{A}^{n}]_{\beta})$ is a filter
basis in the lattice of subspaces of $[\mathcal{A}^{n}]_{\beta}.$
\end{remark}

\begin{proof}
$0\notin\mathfrak{B}_{\beta}$ and for $[M]_{\beta},[N]_{\beta}\in
\mathfrak{B}_{\beta}$, it holds that $[M]_{\beta}\cap\lbrack N]_{\beta
}\supseteq\lbrack M\wedge N]_{\beta}\in\mathfrak{B}_{\beta}$.
\end{proof}

\ 

Let $\widetilde{\mathfrak{B}}_{\beta}$ be a quasipoint of $\mathbb{L}%
([\mathcal{A}^{n}]_{\beta})$ containing $\mathfrak{B}_{\beta}$. There is
exactly one line $[a_{0}]_{\beta}[\mathcal{A}]_{\beta}\in\widetilde
{\mathfrak{B}}_{\beta}$ such that%
\[
\forall M\in\mathfrak{B}:[a_{0}]_{\beta}[\mathcal{A}]_{\beta}\subseteq\lbrack
M]_{\beta},
\]
because the discussion concerning atomic quasipoints for the lattice of
subspaces of finite dimensional vector spaces holds (see Prop.
\ref{PQPsOfTypeI_nFactors}). Of course we have $a_{0}\nsim_{\beta}0$. It
remains to show that $a_{0}\mathcal{A}\wedge M\neq0$ holds for all
$M\in\mathfrak{B}$. We have $P_{a_{0}\mathcal{A}}=E_{\widetilde{a_{0}}}$ and,
since $p_{M}\widetilde{a_{0}}\mathcal{A}\subseteq M$ holds for some
appropriate $p_{M}\in\beta$, we obtain%
\begin{align*}
&  E_{p_{M}\widetilde{a_{0}}}\leq P_{M}\\
\Longleftrightarrow &  p_{M}E_{\widetilde{a_{0}}}\leq P_{M}.
\end{align*}
From this, it follows (since $a_{0}\nsim0$, we have $p_{M}\widetilde{a}%
_{0}\neq0$ for all $p_{M}\in\beta$):%
\[
P_{M}\wedge E_{\widetilde{a_{0}}}\geq p_{M}E_{\widetilde{a_{0}}}\wedge
E_{\widetilde{a_{0}}}=p_{M}E_{\widetilde{a_{0}}}>0,
\]
so $P_{M}\wedge E_{\widetilde{a_{0}}}\neq0$, that is, $a_{0}\mathcal{A}\wedge
M\neq0$. Since $\mathfrak{B}$ is a quasipoint, $a_{0}\mathcal{A}%
\in\mathfrak{B}$ holds from maximality. Summing up, we have shown:

\begin{theorem}
\label{TTypeI_nAllQPab}All quasipoints of $\mathbb{M}_{n}(\mathcal{A)}$ are abelian.
\end{theorem}

The following remark clarifies the relation between several abelian
projections in a single quasipoint:

\begin{remark}
Let $E_{a},E_{b}\in\mathfrak{B}$ be abelian projections. There is some
$r\in\mathcal{P(A)}$ such that $rE_{a}=rE_{b}$.
\end{remark}

\begin{proof}
$E_{c}:=E_{a}\wedge E_{b}$ is an abelian projection in $\mathfrak{B}$, so
there are $p,q\in\mathcal{P(A)}$ such that $E_{c}=pE_{a}=qE_{b}$. Then
$E_{c}=pE_{a}=p^{2}E_{a}=pqE_{b}$ holds and also $E_{c}=pqE_{a}$. If
$p\notin\beta$, then $1-p\in\beta$ and thus%
\[
0=(1-p)pE_{a}=(1-p)E_{c}\in\mathfrak{B},
\]
contradicting $0\notin\mathfrak{B}$. Hence $p\in\beta$ and also $q\in\beta$,
so $pq\in\beta$ and%
\[
pqE_{a}=pqE_{b}.
\]

\end{proof}

\ 

If $\mathcal{R}$ is a type $I_{n}$ algebra with trivial center, i.e.
$\mathcal{R}\simeq\mathbb{M}_{n}(\mathbb{C)\simeq}\mathcal{L(}\mathbb{C}^{n}%
)$, then there are only atomic quasipoints (see Prop.
\ref{PQPsOfTypeI_nFactors} again). An atomic quasipoint is of the form%
\[
\mathfrak{B}_{\mathbb{C}x}:=\{P\in\mathcal{P(R)}\ |\ P_{\mathbb{C}x}\leq P\},
\]
where $x\in\mathcal{H}\backslash\{0\}$. Of course, $P_{\mathbb{C}x}$ is an
abelian projection, too. The Stone spectrum $\mathcal{Q(R)}$ is discrete in
this case, since atomic quasipoints are isolated points of the Stone spectrum
$\mathcal{Q(R)}$.

\ 

As a corollary of the results proved above, we obtain the main result:

\begin{theorem}
\label{TQPsOfCenterParamUnitaryOrbits}Let $\mathcal{R}$ be a type $I_{n}$ von
Neumann algebra. The quasipoints of $\mathcal{A}=\mathcal{C(R)}$ parametrize
the orbits of the unitary group $\mathcal{U(R)}$ acting on the Stone spectrum
$\mathcal{Q(R)}$ of $\mathcal{R}$.
\end{theorem}

\begin{proof}
All quasipoints of $\mathcal{R}$ are abelian (Thm. \ref{TTypeI_nAllQPab}),
i.e. we have $\mathcal{Q(R)}=\mathcal{Q}^{ab}(\mathcal{A})$. Using this, we
apply Thm. \ref{TAbQPcQP}: the mapping%
\begin{align*}
\zeta:\mathcal{Q}^{ab}\mathcal{(R)}  &  \longrightarrow\mathcal{Q(A)},\\
\mathfrak{B}  &  \longmapsto\mathfrak{B}\cap\mathcal{A}%
\end{align*}
is surjective. Since $\mathcal{R}$ is finite, we can replace partial
isometries by unitary operators in Thm. \ref{TAbQPcQP}. Therefore,
$\zeta(\mathfrak{B})=\zeta(\widetilde{\mathfrak{B}})$ holds if and only if
there is a \emph{unitary} $U\in\mathcal{R}$ such that $U.\mathfrak{B}%
=\widetilde{\mathfrak{B}}$ (see Def. \ref{DUActionOnQP}). It follows that the
quasipoints of $\mathcal{A}=\mathcal{C(R)}$ parametrize the orbits of the
unitary group $\mathcal{U(R)}$ acting on $\mathcal{Q(R)}$.
\end{proof}

\section{Acknowledgement}

I would like to thank H. de Groote for many helpful discussions and
substantial support. Financial support by the Studienstiftung des Deutschen
Volkes is gratefully acknowledged.

\end{document}